\documentclass[11pt]{amsart}
\usepackage{amssymb}
\usepackage{enumitem}   

\theoremstyle{cupthm}
\newtheorem{theorem}{Theorem}[section]
\newtheorem{proposition}[theorem]{Proposition}
\newtheorem{lemma}[theorem]{Lemma}
\newtheorem{corollary}[theorem]{Corollary}

\newtheorem{definition}[theorem]{Definition}
\newtheorem{example}[theorem]{Example}
\newtheorem{remark}[theorem]{Remark}

\newtheorem*{notation}{Notation}

\def\c{Con(\frak g,\mathcal E)}

\def\ge{(\frak g,\mathcal E)}

      \newcommand{\nil}{\emptyset}
      \newcommand{\g}{{\mathfrak g}}
      \newcommand{\h}{{\mathfrak{h}}}

       \newcommand{\E}{{\mathcal{E}}}
       \newcommand{\F}{{\mathcal{F}}}
      
     \newcommand{\Con}{\operatorname{Con}}
     
      \newcommand{\con}[2]{\Con(\mathfrak{#1},\mathcal{#2})}
       
      \newcommand{\tcon}[2]{\Con_{\text{t}}(\mathfrak{#1},\mathcal{#2})}

\begin{document}
\title{A survey of configurations }

\author[Akram Yousofzadeh]{ Akram Yousofzadeh}


\begin{abstract}
To any finite ordered subset and any finite partition of a group a set of tuples of positive integers, named as configurations, is associated that  describes the group's behavior. The present paper provides an exposition
 of this 
notion and related topics, including amenability, paradoxical decomposition and the group theoretical properties.
\end{abstract}



\maketitle

\section{Introduction and preliminaries}

The theory of configurations that we are going to describe, began in 2001 when J. M. Rosenblatt G. A. Willis  published a paper entitled "Weak
Convergence is Not Strong Convergence for Amenable groups". Using some features of the Cayley graph of a group, coloring this graph and giving an equivalence condition for amenability, they defined the notion of configuration and presented the first construction of nets which are weak-* converging to invariance but not strongly converging to invariance. 
 This concept can be extended to other algebraic structures such as semigroups and hypergroups. It has results in amenability of a class of Banach algebras, as well.\\

We start with a brief overview of the theory, by giving the definitions and fundamental theorems. \\

Let $G$ be a group, $\frak g=\{g_1,\dots, g_n\}$ be a finite subset and $\mathcal E=\{E_1,\dots,E_m\}$ be a finite partition of $G$. The oriented Cayley graph $Cay(G,\frak g)$ corresponding to $\frak g$, has vertices being the elements of $G$ and directed edges connecting the element $h$ to the element $g_ih$ for some $i\in\{1,\dots,n\}$. Color the vertices of $Cay(G,\frak g)$, such that all vertices in $E_i$ are given the same distinct color. This way, the Cayley graph is categorized into finitely many classes which define  {\it{configurations}}.\\

Each
configuration corresponding to the finite string $\frak{g}=(g_1,\dots, g_n)$ and the finite partition $\mathcal E=\{E_1,\dots,E_m\}$ of $G$ is an
$(n+1)$-tuple $C=(C_0,C_1,\dots, C_n)$ with $C_k\in \{1,\dots,m\}$
such that there is a set of elements $x_0,x_1\dots x_n\in G$ which $x_0\in E_{C_0}$ and 
for any $j=1,\dots,n$ we have $x_j=g_jx_0\in E_{C_j}.$  The element $x_0$ is
called a base point  and each $x_j$ is called a branch point for
$C$. Obviously the
elements $x_0,x_1\dots x_n$ are not necessarily unique. Despite the fact that there may be infinitely many tuples $(x_0,x_1,\dots,x_{n})$ which define a certain configuration, there are only a finite number of configurations, for any choice of $\frak g$ and $\mathcal E$.\\

Let $C$ be a configuration. Consider the set of all base points or equivalently the set of all elements of $G$ that can be placed in the position of $x_0$ and denote this set by $x_0(C)$. It is easy to see that different configurations have disjoint sets of base points. These sets are the key points to almost all proofs.
\begin{example}
Let $G=\mathbb{Z}$ be the additive group of integer numbers. Put $g_1=1, g_2=2$ and $\mathcal E=\{E_1, E_2,E_3\},$ where $E_1=3\mathbb{Z}$, $E_2=3\mathbb{Z}+1$, and $E_3=3\mathbb{Z}+2$. Then the set of configurations corresponding to $\frak g$ and $\mathcal E$ consists of
$$C_1=(1,2,3), C_2=(2,3,1)\ \ \text{and}\ C_3=(3,1,2).$$
Here we have $x_0(C_1)=E_1, x_0(C_2)=E_2$ and $x_0(C_3)=E_3.$
\end{example}

 Let $P(G)$ denote the set of positive functions $f \in L^1(G)$ such that $\int_G f d\lambda = 1$. If $f\in P(G)$, define $f_C$ to be the sum of all values of $f$ on the set $x_0(C)$, i.e., $$f_C=\sum\{f(x),\ x\in x_0(c)\}.$$ \\
 
 With the configurations corresponding to the pair $\ge$ there is associated a system of equations, stating that for the color $i$ and for each $j\in\{1,\dots,n\}$, the sum of all $f_C$ with base point of $i$-th color, should be the sum of all  $f_C$ with $j$-th branch point and the same color $i$. This system is called the system of configuration equations and is denoted by $Eq\ge$. Indeed this system is defined through
$$\sum\{f_C,\ x_0(C)\subseteq E_i\}=\sum\{f_C,\ x_j(C)\subseteq E_i\},\ \ \ \ \ \ (1\leq j\leq n,\ 1\leq i\leq m).$$
It is proved that the solutions to this system depend only on the values of the string $(f_C)$ and not on the function $f\in P(G).$ Therefor to solve the system we shall only look for strings $(f_C)$ of real numbers, satisfying the equations. This solution is called normalized if $f_C\geq 0$ for each $C$ and $\sum f_C=1.$\\

A locally compact group $G$ is called left amenable if it admits a left
invariant mean on $L^{\infty}(G).$ That is a positive linear functional   $M$ on  $L^{\infty}(G)$ such that
$M(1)=1$, which is invariant; i.e.,   $$M(_af)=M(f),\ \ \ \ \ (a\in G, f\in L^{\infty}(G)).$$\\

The amenability of groups  is intimately close to the existence of normalized solutions for configuration equations. Indeed we have the following theorem which is the essence of many results proved during the years after 2001.\\

 \begin{theorem}\cite[Theorem 2.4]{rw} \label{main}
 There is a normalized solution of every possible instance of the configuration
equations if and only if $G$ is amenable.
\end{theorem}

 Clearly, the word normalized can be replaced by nonzero nonnegative, here. But even more has been proved: There is a nonzero solution of every possible instance of the configuration
equations if and only if $G$ is amenable (see \cite{ar}).\\

The aim of the work of Rosenblatt and Willis was to present appropriate nets
which are weak* converging to invariance; i.e., a net $(f_{\alpha})$ in $P(G)$ such that $\langle f_{\alpha}-_gf_{\alpha}, F\rangle$ tends to zero for any $F\in L^{\infty}(G)$ and every $g\in G$.
If $G$ is an infinite locally compact group, varying partitions $\mathcal E$ and finite subsets $\frak g$ of $G$, and taking normalized solutions to each $Eq(\frak g,\mathcal E),$ the desired net $(f_{\alpha})$ in $P(G)$ can be constructed. The statement of the theorem is

\begin{theorem} \cite[Theorem 2.9]{rw}
If $G$ is an infinite discrete amenable group, then there exists a net $(f_{\alpha})$
in $P(G)$ converging weak* to invariance,
such that for any $e\neq x\in G,$  $\|f_{\alpha}-_xf_{\alpha}\|_1=2$ .
\end{theorem}

Systems with no normalized solution provide us a useful tool to investigate and compare non-amenable groups. We shall give an explanation in section 4 below.\\

The application  of the theory of configuration to group theory is given in sections  2 and 3.
In section 4 we apply configurations to construct paradoxical decompositions of non-amenable groups and to determine new upper bounds for their Tarski's numbers. Section 5 is devoted to a short  summary of studies about the configuration of semigroups and hypergroups as well as Lau algebras.

\section{ group properties  characterized by configurations}
Amenability of a discrete group $G$ is equivalent to the amenability of all finitely generated subgroups of $G$. Throughout this section we assume that the group $G$ is finitely generated and we state the properties that can be characterized through the coincidence of the set of  configurations.\\

 If $\frak g$ and $\mathcal E$ are as in the pervious section and $\frak g$ generates the group $G$, we call $\ge$ a configuration pair. The set of all configurations corresponding to $\ge$ is denoted by $\c$ and is called a configuration set.  By $Con(G)$ we mean the family of all possible configuration sets of $G$.\\

It is useful sometimes to go from one element of $Con(G)$ to another one, by switching the generator $\frak g$ or the partition $\mathcal E$.  For example it is important in many considerations below, that a special element of $G$ (mostly $e_G$) belongs to a specified member of $\mathcal E$ so that the configuration set $\c$ does not differ. It is actually possible.
Indeed, for any partition $\mathcal E$ and element $x\in G$ we have $Con(\frak g, \mathcal E)=Con(\frak g, \mathcal E')$, where $\mathcal E'=\{E_1x,\dots,E_mx\}$. It is also proved that if $\mathcal F$
is refinement of the partition $\mathcal E$ and $Eq(\frak g,\mathcal F)$ admits a normalized solution, then so is  $Eq(\frak g,\mathcal E)$, see \cite{farsi}.\\

We say that two groups $G$ and $H$ are configuration equivalent and write $G\approx H,$ if
$$Con(G)=Con(H).$$
 Several properties of groups is clearly passed through the configuration equivalence. For instance, configuration equivalence groups have the same minimum number of generators. Also according to the Theorem \ref{main}, one of them is amenable if and only if the other one is amenable.\\

In the following cases, an appropriate selection of a configuration set $\c$ would show that two specified configuration equivalent groups are isomorphic.\\

Let $G$ be finite. Selecting $\frak g$ the ordered generating set of $G$ and $\mathcal E$ the most refined partition of $G$, i.e., the set of all singletons of elements of $G$, an elementary proof shows that $G\approx H$ implies $G\cong H$ (see the proof of \cite[proposition 4.3]{arw}).
Configuration equivalence can also specify some other group theoretical properties of groups. We list them as follow referring directly to \cite{arw} and \cite{trya}.

\begin{theorem}
Let $\mathbb{F}_n$ be the free group of positive rank $n$. If $H$ is a
finitely generated group such that $H\approx \mathbb{F}_n$, then $H\cong \mathbb{F}_n$.
\end{theorem}

If  for a periodic group $H$,  $G\approx H,$  then $G$ must be periodic as well. Indeed
\begin{theorem}
Let $G$ be a finitely generated group having an element
of infinite order. Then there is a partition, $\mathcal E$, of $G$ and a generating set
$\frak g=(g_1,\dots, g_n)$ such that $\c$ is not a configuration set of any periodic
group.
\end{theorem}

Let $n$ be a positive integer and  $S$ be the free semigroup
on the set $\{x_1,\dots,x_n\}$. Suppose that
$\mu= \mu(x_1,\dots,x_n)$ and $\nu = \nu(x_1,\dots,x_n)$ are two elements in $S$. We say that
$\mu =\nu$ is a semigroup law in a group $G$ if for any $n$-tuple $(g_1,\dots,g_n)$ of
elements of $G$, we have $\mu(g_1,\dots,g_n) =\nu(g_1,\dots,g_n).$ For example $G$ is abelian if and only if the semigroup law $xy=yx$ is satisfied.
The following theorem implies the stability of certain group properties for configuration equivalent groups.\\
\begin{theorem}
Let $G_1$ and $G_2$ be two finitely generated groups with $G_1\approx G_2$. Then  $G_1$ and $G_2$ satisfy the same semigroup laws.
\end{theorem}

So if $G_1$ is abelian, then so is $G_2$.  For this case, the isomorphism is also concluded:\\
\begin{theorem}
Let $G$ be a finitely generated non-abelian group. Then
there is a partition, $\mathcal E$, of $G$ and a generating set $\{g_1, g_2,\dots, g_n\}$ such that the
corresponding configuration set cannot arise from an abelian group.
Moreover if $G_1$ and $G_2$ are finitely generated abelian groups such
that $G_1\approx G_2$, then $G_1\cong G_2$.
\end{theorem}
As the nilpotency of groups of a certain class can be stated by the language of semigroup laws (see \cite{taylor}), we have the following
\begin{theorem}\label{nil}
Let $G_1$ and $G_2$ be two finitely generated groups with $G_1\approx G_2$. Suppose that $G_1$ is nilpotent of class $c$. Then so is $G_2$.
\end{theorem}
This result has been improved recently. It is proved in \cite{rs2} that two configuration equivalent groups obey the same  {group laws}.\\

\begin{proposition}
Let $G$ and $H$ be two groups with  $G\approx H$ and suppose that $G$ satisfies the {group law} $\mu(x_1, . . . , x_n) =e_G$. Then $H$ satisfies the same law.
\end{proposition}

\begin{corollary}
Let $G$ and $H$ be two groups such that $G\approx H$. Then the solubility of $G$ is equivalent to $H$'s.
Furthermore, their derived lengths are the same.
\end{corollary}

For the rest of this section, we mention the results concerning the quotient groups. For a group $G$, we denote by $\mathcal F(G)$ and $\mathcal A(G)$ the set of isomorphism classes of finite and abelian quotients of $G$, respectively. Let us begin with the following two theorems from \cite{arw} and \cite{ary}.

\begin{theorem}\label{fq}
Let $G_1$ and $G_2$ be two finitely generated groups with $G_1\approx G_2$. Then $G_1$ and $G_2$ have the same set of finite quotients. In other words, if $G_1$ has a normal subgroup, $N_1$, with finite index. Then $G_2$
has a normal subgroup $N_2$ with $G_2/N_2\cong G_1/N_1$; i.e., $\mathcal F(G_1)=\mathcal F(G_2).$
\end{theorem}

\begin{proposition}\label{abq} Let $G$ and $H$ be two finitely generated groups such that $G\approx H$.
Then $\mathcal A(G_1)=\mathcal A(G_2)$.
\end{proposition}

As a direct application of Theorem \ref{abq}, we point out that for two arbitrary configuration equivalent groups $G$ and $H$, $\frac{G}{G'}\cong\frac{H}{H'}$ (See \cite{trya}).\\

Let $G$ be a finitely generated torsion-free nilpotent group of class 2. Denote
by $T=T(G)$ the isolator of $G'$ in $G$, i.e., the set of all elements $x\in G$ such that
$x^s\in G'$ for some nonzero integer $s$. The set $T$ is a central subgroup of $G$ and $T/G'$
is the torsion subgroup of $G/G'$. Then there exist positive integers $m$ and $n$ such
that $G/T\cong \mathbb{Z}^n$ and $T\cong \mathbb{Z}^m$. The class of all such groups $G$ are denoted by
$\mathcal I(n,m)$ (see \cite{DSegal}, p. 260). If $G\in \mathcal I(n,m)$ and $G\approx H$, then $H$ is a finitely generated torsion-free nilpotent group of class 2, as well. So $H\in \mathcal I(n',m')$ for some positive integers $n'$ and $m'$. It is seen that $n=n'$ and $m=m'.$
For the particular case, when $(n,m) \in \{(3, 2), (3, 3), (n, 1)\},$ it is well known that if $G, H\in \mathcal I(n,m)$ and $G$ and $H$ are not isomorphic, then they have different classes of finite quotients. Thus the following theorem obtained.

\begin{theorem} Let $G\in \mathcal I(n,m)$ and $G\approx H$ for some finitely generated group $H$.
If $(n,m) \in \{(3, 2), (3, 3), (n, 1)\},$ then $G\cong H.$
\end{theorem}

Note that $Tr(3, \mathbb{Z})\in \mathcal I(2,1)$. So if $G$ is a finitely generated group such that $G \approx Tr(3, \mathbb{Z})$, then $G\cong Tr(3, \mathbb{Z}).$\\

Following  Theorem \ref{nil}, it is naturally asked wether the configuration equivalence of two nilpotent groups necessities them to be isomorphic. This question is investigated for nilpotent $FC$-groups in \cite{ary}.\\

An $FC$-group is a group in which each element has finitely many
conjugates. A finitely generated group $G$ is an $FC$-group if and only if
$G/Z(G)$ is finite, where $Z(G)$ is the center of $G$. Also If $G$ is a finitely generated $FC$-group,
then $G$ is isomorphic with a subgroup of $\mathbb{Z}n \times F$, for some finite group $F$. (see \cite{tomkinson}).
It is shown in \cite{tomkinson} that for each finitely generated $FC$-group $G$, the torsion subgroup
$Tor(G)$, and commutator subgroup $G'$ of $G$ are finite and therefore the factor group $G/Tor(G)$ is an abelian torsion free group.\\

For a group $G$ the non cancelation set, $\chi(G)$ is the set of
isomorphism classes of groups $H$ such that $G\times \mathbb Z\cong H\times \mathbb Z.$ If the commutator subgroups of the groups $G$ and $H$
are finite, then $\mathcal F(G)=\mathcal F(H)$ necessities them to have the same cancelation sets (see \cite[Theorem 2.1]{witbooi}). So configuration equivalent groups have equal cancelation sets. Using these facts the following theorem is followed\\

\begin{theorem}
Let $Con(G)=Con(\mathbb{Z}_n \times F)$, where $F$ is a finite group and $G$ is
an $FC$-group. Then $G\cong \mathbb{Z}_n \times F$.
\end{theorem}

 A polycyclic group is a solvable group that satisfies the maximal condition on subgroups (that is, every subgroup is finitely generated). Polycyclic groups are finitely presented, and this makes them interesting to work. Equivalently, a group $G$ is polycyclic if and only if it admits a subnormal series with cyclic factors, that is a finite set of subgroups, say $G_0,\dots,G_n,$ such that
$G_0$ coincides with $G$, $G_n$ is the trivial subgroup, $G_{i+1}$ is a normal subgroup of $G_i$ (for every $i$ between 0 and $n-1$) and the quotient group $G_i/G_{i+1}$ is a cyclic group.\\

\begin{proposition} Let $G$ and $H$ be two polycyclic groups such that $G\approx H.$ If $G$
admits a normal series with infinite cyclic factors, then $H$ is torsion free.
\end{proposition}

In a { {polycyclic}} group $G$ the number of infinite factors in a cyclic series, which
is known as the {\it {Hirsch}} length, is independent of the choice of series.
By induction on the nilpotency class of a torsion free nilpotent group $G$, one can see that if $G\approx H$, then
$H$ is also a torsion free nilpotent group with the same Hirsch length of $G$:\\

\begin{theorem}\label{nilfc}
Let $G$ be a finitely generated nilpotent $FC$-group and $G\approx H$. Then
$H$ is a nilpotent $FC$-group too. Furthermore $Z(G)\cong Z(H)$ and $G/Z(G)\cong H/Z(H)$. Also
$G$ and $H$ have the same Hircsh lengths.
\end{theorem}

\begin{definition}\em
Let $\mathbb F_m$ denote the free group of rank $m$ on set $\{x_1,\dots,x_m\}.$ For a fixed $n$,
let $\mathbb F^n_m$ denote the subgroup of $\mathbb F_m$ generated by $g^n$ for all $g\in \mathbb F_m$. Then $\mathbb F^n_m$ is a normal
subgroup of $\mathbb F_m$ and we define the Burnside group $B(m, n)$ to be the quotient group $\frac{\mathbb F_m}{\mathbb F^n_m}$.
\end{definition}
Burnside groups can also presented via generators and relations.  Novikov and Adyan presented the group $B(m,n)$ as follows
(see \cite{robin})
$$B(m,n)=\langle x_1,\dots, x_m; w^n = 1 \rangle.$$ Using this form, we have the next theorem

\begin{theorem}
The following statements are equivalent:
\begin{enumerate}
\item Con(B(m, n)) = Con(B(m', n')).
 \item m = m', n = n'.
\item B(m, n) = B(m', n').
\end{enumerate}
\end{theorem}

\vspace{.3cm}

The concept of golden configuration pair is defined through the next notation and definition which is initiated by the authors in \cite{rs2}.
The existence of a golden pair can guaranty the isomorphism of some sorts of groups, provided that they have the same set of configurations (see proposition \ref{last} below).

\begin{notation}
  \label{not}
Let $G$ be a group with $\mathfrak{g}=(g_1,\dotsc,g_n)$ as its generating set. Let $p\in\mathbb{N}$, $J\in\{1,2,\dotsc,n\}^p$ and $\sigma\in\{\pm1\}^p$. Denote the product
$\prod_{i=1}^Ng_{J(i)}^{\sigma(i)}$ by $W(J,\sigma;\mathfrak{g})$. The pair
$(J,\sigma)$ is called a \textit{representative pair} on
$\g$ and $W(J,\sigma;\mathfrak{g})$ a \textit{word} corresponding to $(J,\sigma)$ in $\g$. The representative pair $(J,\sigma)$ is called \textit{reduced}, if $\sigma(i)=\sigma(i+1)$, whenever
$J(i)=J(i+1)$.
\end{notation}

\begin{definition}\em
\label{golden}
Let $G$ be a group. We call a configuration pair $(\g,\E)$ golden, if for each configuration
 pair $(\g',\E')$ for $G$ with $\E'=\{E_1',\dotsc,E'_m\}$, $e_G\in E_1'$ and $\con gE=\con {g'}{E'}$, we have the following implication:
$$ W(J,\sigma;\g)\neq e_G \Rightarrow W (J,\sigma;\g')\not\in E'_1$$
in which $(J,\sigma)$ is a representative pair.
\end{definition}
\noindent For the details of the proof of the following proposition see \cite{rs2}.
\begin{proposition}
 \label{last}
  Let $G$ be a Hopfian group with a golden configuration pair and $H$ be a group such that $G\approx H$. Then $G$ is finitely presented if and only if $H$ is finitely presented, and in the case that $G$ has a finite set of defining relators, we have $G\cong H$.
  \end{proposition}

Two groups $G$ and $H$ are \textit{strong configuration equivalence} written as $(G;\mathfrak{g})\approx_s(H;\mathfrak{h})$, if there exist ordered generating sets $\mathfrak{g}$ and
 $\mathfrak{h}$ for $G$ and $H$ respectively, such that

\noindent 1) For each partition $\mathcal{E}$ of $G$ there exists  a partition $\mathcal{F}$ of $H$ such that $\Con(\mathfrak{g},\mathcal{E})=\Con(\mathfrak{h},\mathcal{F})$ and
 2) For each partition  $\mathcal{F}$ of $H$ there is a partition $\mathcal{E}$ of $G$ such that
 $\Con(\mathfrak{h},\mathcal{F})=\Con(\mathfrak{g},\mathcal{E})$.

In \cite{rs2} it is proved that this version of configuration equivalence for two groups implies the isomorphism of them:
 \begin{theorem}\cite[Theorem 4.3]{rs2}
 Two groups are strong configuration equivalence if and only if they are isomorphic.
 \end{theorem}


\section{Two-sided configurations}

To investigate the behavior of a group in some cases, we need to know how a special element acts from both left and right sides. In this
section, we give a new type of configurations which depends on the left and right
translations and we shall give a review of results  from \cite{twos}.

\begin{definition} \em Let $G$ be a finitely generated group. Given a generating sequence
$\mathfrak g = (g_1,\dots, g_n)$ and a partition $\mathcal E =\{E_1,\dots,E_m\}$ of $G$, a two-sided configuration
corresponding to the pair $\ge$ is an $(2n + 1)$-tuple $C=(C_0, C_1,\dots, C_{2n})$ satisfying the following two conditions\\
\begin{enumerate}
\item $C_i\in\{1,\dots,m\}$ for each $0\leq i\leq 2n$, and
 \item there exists $x\in E_{C_0}$ 
 such that
$g_ix\in E_{C_i}$ and $xg_i\in E_{C_{i+n}}$ for each $i \in \{1,\dots, n\}$.
 \end{enumerate}
\end{definition}

We denote by $Con'(\frak g, \mathcal E)$ the set of all two-sided configurations corresponding to
$\ge$, and let $$Con'(G) = \{Con'(\frak g, \mathcal E) :\ \ \ge \text{\ is\ a\ configuration\ pair\ for\ G}\}.$$\\

It is worth pointing out that two definitions of configurations are different (see \cite{twos}). In fact there is no direct correspondence between one-sided and two-sided
configuration sets for the given generator and partition. It is naturally asked
if there exist groups $G_1$ and $G_2$ so that $G_1\approx G_2$ but $Con'(G_1)\neq Con'(G_2).$\\

Two groups $G_1$ and $G_2$ are called two-sided configuration equivalent if $Con'(G_1)=Con'(G_2).$ It is obvious that if $Con'(G_1) = Con'(G_2)$, then $Con(G_1) = Con(G_2)$.
Therefore, every conclusion obtained for
configuration equivalent groups is true for two-sided configuration equivalent groups, as well.
Clearly if $G_1$ and $G_2$ are abelian,  then the two concepts are the same.\\

In the pervious section we saw that the equivalence of configurations of two finitely generated groups, implies that they have the same finite quotients. It is not difficult to find examples to prove that the converse is not true. In the following, two polycyclic groups are given in  this direction.\\

Let $K=\mathbb Q(\sqrt {10})$ be as a number field of degree 2, and $O =\mathbb Z +\sqrt{10}\mathbb Z$ the ring
of integers of $K$. It is to be noted that  $O_1 =3 \mathbb Z+(1+\sqrt{10})\mathbb Z$ and $O_2 = O$ are ideals which are in
distinct ideal classes of $O$. The element $u = 3+\sqrt10$  is a unit element of infinite order
in $O$. For $i = 1, 2,$ choose a $\mathbb Z$-basis for $O_i$, and let $x_i$ be the matrix corresponding
to multiplication by $u$ with respect to this basis. Accordingly
\begin{center}
$x_1 =\left[ \begin{array}{cc}
2 &3\\
3 & 4
\end{array} \right]$ and $ \ \ \ x_1 =\left[ \begin{array}{cc}
3 &10\\
1 & 3
\end{array} \right].$
\end{center}
 As members of  $GL2(Z),$ neither $x_1$ are $x_2$ nor $x_1-1$ and $x_2$ are similar. So $\langle x_1 \rangle$ and $\langle x_2 \rangle$ are not conjugate in $GL_2(Z).$ For details,
see \cite[p. 256]{DSegal}.

\begin{lemma} Let $G_i=\mathbb Z_2\otimes_{\rho_i}\langle x_i\rangle$, where $\rho_i\in Aut(\mathbb Z_2)$ and $\rho_i (m, n) = x_i(m, n)$
for $i = 1, 2$ and $m, n \in\mathbb Z$. Then $G_1$ is not isomorphic to $G_2$, but $\mathcal F(G_1) =\mathcal F(G_2).$
\end{lemma}
Note that these groups are torsion free polycyclic. 
This observation leads us to have the following theorem.

\begin{theorem} \cite{twos} There exist non-isomorphic polycyclic torsion free groups $G_1$ and
$G_2$ such that $Con'(G_1)\neq Con'(G_2)$ and $\mathcal F(G_1) =\mathcal F(G_2)$.
\end{theorem}

Let $G_1* G_2$ denote the free product of two groups $G_1$ and $G_2$.
Each element $g$ of $G_1*G_2$ can be uniquely written in the form $g = g_1g_2\dots g_r$, where
$r\geq 0,\ 1 \neq g_i\in G_{n_i},$ $n_i\in \{1, 2\}$ and $n_i \neq n_{i+1}$. The case $r = 0$ is interpreted as
$g = 1$ (see \cite[pp. 161-178]{robin}).\\

Let $G = F_m$ and $H = F_n$ be two finite rank free groups. Then
their free product is again a finite rank free group $F_{m+n}$ (see \cite[p. 164]{robin}), and so if
$Con(K) = Con(G * H)$, then $K\cong  G * H$ by [2].
The following theorem is obtained by direct uses of definition (see \cite{twos}).\\

\begin{theorem}
 Let $F$ be a finite group and $H$ be a finite group or an infinite cyclic group such that $Con'(G) = Con'(H* F)$. Then $G \cong H * F$.
\end{theorem}

In section 2 we recalled some results concerning nilpotent groups. It is well-known that
each finitely generated nilpotent group is isomorphic to a
subgroup of $Tr(n;Z)\times F$ for some positive integer n and a finite nilpotent group
$F$, where $Tr(n;Z)$ is the group of $n\times n$ upper triangular matrices over $\mathbb Z$ (see \cite[p. 158]{rw}). The following theorem, solves the problem of two-sided configuration
equivalence and isomorphisms for this case.

\begin{theorem} \cite[Theorem 3.5]{twos} Let $G$ be a finitely generated group such that
$$Con'(G) = Con'(Tr(n;\mathbb Z) \times F),$$
where $F$ is a finite group. Then $G \cong Tr(n;\mathbb Z) \times F.$
\end{theorem}


Let $N$ be a normal subgroup of a group $G$. A partition $\E$ of $G$, is called a \textit{quotient partition}, if it is an inverse image of a partition of $G/N$ under the quotient map; In other words, there exits a partition $\mathcal P$ of $G/N$, such that {\small $\E=q^{-1}(\mathcal P):=\{q^{-1}(P):\,P\in\mathcal P\}$.} If a refinement $\E'$ of a quotient partition $\E$ of $G$ is itself quotient, we say $\E'$ is a \textit{quotient refinement }of $E$. By a \textit{quotient configuration pair}, we mean a configuration pair $(\hat\g,\E)$ in which $\hat\g$ is a quotient extension generating set and $\E$ is a quotient partition of $G$.\\

\begin{definition} \cite{rs3}\em
\label{recognizable}
Let $G$ be a group with a normal subgroup $N$ and a quotient configuration pair $(\hat\g,\E)$. We say that $(\hat\g,\E)$ is \textit{recognizable} w.r.t $N$, if $\tcon {\hat g}E=\tcon {\hat h}F$, for a configuration pair $(\hat\h,\F)$ of a groups $H$ , then for every $g\in G\setminus N$, there is a representative pair, $(J_g,\sigma_g)$ such that
$$g=W(J_g,\sigma_g;\hat\g)\quad\text{and}\quad W(J_g,\sigma_g;\hat\h)F\cap F=\nil$$
 where $F\in\F$ is in the same color as an element of $\E$ which contains $N$.
 \end{definition}

It is shown that if $G$ is a finitely presented Hopfian group with a recognizable configuration pair, and $H$ be a group with $Con'(G)=Con'(H)$, then $G\cong H$. Also it is proved that every polycyclic group has a recognizable configuration pair. Therefore we have the following important proposition:
\begin{proposition}
\label{poly iso}
Let $G$ be a polycyclic group and  $Con'(G)=Con'(H)$ for a finitely generated group $H$. Then $G\cong H$.
\end{proposition}

 The important point to note here is that every finitely generated nilpotent group is polycyclic Therefore
 for two-sided configuration equivalent groups  $G$ and $H$ such that $G$ is nilpotent, we have $G\cong H$. However, a more general
 result has been obtained:

\begin{theorem} \cite[Theorem 5.3]{rs3}
Let $G$ be a polycyclic--by--finite group, such that $Con'(G)=Con'(H)$ for a finitely generated group $H$. Then $G\cong H$.
\end{theorem}

We investigated the configuration equivalence of nilpotent FC-groups in section 3. These groups have finite commutators. Theorem \ref{nilfc} is extended in the following sense.\\

\begin{theorem}
Let $G$ be a finitely generated group having finite commutator subgroup. Assume that $Con'(G)=Con'(H)$ for a group $H$. Then $G\cong H$.
\end{theorem}
\noindent In particular if $G$ is a finitely generated FC-group, which $Con'(G)=Con'(H)$ for a group $H$. Then $G\cong H$ (see \cite{rs3} and \cite{tomkinson}).\\


The main question of section 2 is
`{ Does the configuration equivalence imply the isomorphism of groups}'.
Recently, this question is answered negatively (see \cite{rs1}). In that paper the authors give an example including two soluble groups which are two-sided configuration equivalent (and hence configuration equivalent) but they are not isomorphic. Therefor the main question is not correct even for amenable groups.

\begin{theorem}\label{non-isomorphic generatin-bijective groups}
There exist non-isomorphic finitely generating groups with the same two-sided configuration sets.
\end{theorem}

\section{paradoxical decompositions}
 Let $G$ be a discrete group. The group $G$ is amenable if it admits a finitely
additive probability measure $\mu: \mathcal P(G)\rightarrow [0,1]$ on the power set of $G$, which is left invariant. 
In other words, for any $g\in G$ and any $A\subseteq G$, $\mu(gA) = \mu(A),$
where $gA = \{ga:\ a\in A\}.$ 

\begin{definition}\em
 Let $G$ be a discrete group. $G$ admits a paradoxical decomposition if there exist $m,n\in \mathbb N,$ disjoint subsets $A_1,\dots,A_n$, $B_1,\dots,B_m$ of $G$ and subsets $\{g_1,\dots,g_n\}$ and $h_1,\dots,h_m$ of $G$ such that
 $$G=\bigcup_{i=1}^ng_iA_i=\bigcup_{j=1}^mh_jB_j.$$
\end{definition}

The number $\tau = n + m$ for $n$ and $m$ in the above definition is called the
Tarski number of that paradoxical decomposition; the minimum of all such numbers
over all the possible paradoxical decompositions of $G$, is called the Tarski
number of $G$ and is denoted by $\tau(G)$. We put $\tau(G) =\infty$ in the case that there is no paradoxical decomposition, (see \cite{robin}).
It is obvious that $m,n\geq 2.$ Hence for any group $G$ we have $\tau(G)\geq 4$. The groups with Tarski number 4 are completely determined; in fact $\tau(G)=4$ if and only if $G$ contains a non-abelian free subgroup \cite[5.8.38]{sapir}.
The first known groups with Tarski numbers 5 and 6 are given in \cite{sapir ershov}. It is to be noted that no group with finite exact Tarski number $\geq 7$ has been presented so far.  However, the set of all finite Tarski numbers is infinite (see \cite{sapir ershov}).  \\

It is easily seen that for two configuration equivalent groups $G_1$ and $G_2$,  $G_1$ is amenable if and only
 if $G_2$ is amenable. In other words, $\tau(G_1)=\infty$ if and only if $\tau(G_2)=\infty$. The next theorem states that this
  result is true even for non-amenable groups. 

\begin{theorem} \cite[Theorem 3.8]{ytr} Let $G_1$ and $G_2$ be two groups such that $Con(G_1)=Con(G_2).$ Then
$\tau(G_1)=\tau(G_2).$
\end{theorem}
This concludes in particular that one of the configuration equivalent groups contains a non-abelian free group if and only
if the other one does (see, for example, \cite[Theorem 5.8.38]{sapir}). For subgroups of two configuration equivalent groups, we have the following theorem \cite{farsi}.
\begin{theorem}
Let $G_1\approx G_2$ and $G_1$ has a normal subgroup $H_1$ of index $p<\infty.$ Then there is a normal subgroup $H_2$ of $G_2$ such that 
$$\frac{1}{p}\leq \frac{\tau(H_1)-2}{\tau(H_2)-2}\leq p.$$\\
For example if $p=2$ and $\tau(H_1)=5,$ then $4\leq\tau(H_2)\leq 8$.
\end{theorem}

It is clear that if $G$ admits a paradoxical decomposition, then it is not amenable. The Tarski alternative states that the converse is also true \cite{wagon}. Now the problem is 'what is the relation between paradoxical decompositions of a group and its configurations'.\\

In the rest of this section we recall two procedures that construct paradoxical decompositions for a given group, using a
system of configuration equations with no nonnegative nonzero solution. We omit the proofs but for the sake of importance,
we point out that the proofs are strongly based on the following two lemmas. For details see \cite{y-r}.\\

Let $(\frak g, \mathcal E)$ is fixed and put
$$A^j_i = \{C\in Con(g, E) :\ x_j(C) \subseteq Ei\},\ \ \  (1\leq i\leq m,\ 1\leq j\leq n).$$
It is clear that for each $j\in \{0,\dots, n\}$, ${A^j_1,\dots, A^j_m }$ is a partition of $Con(g, E).$\\

\begin{lemma} For each $i\in \{1,\dots,m\}$ and $j,j'\in \{0,\dots, n\}$,
$$g_{j'}^{-1} g_j(\bigcup_{C\in A^j_i}x_0(C))=(\bigcup_{C\in A^{j'}_i}x_0(C)).$$
\end{lemma}

\begin{lemma} \label{ll} If the system $AX=0$ has no nonnegative nonzero solutions, then by row operations, $A$ can be changed
into an equivalent matrix $B$ with nonnegative entries and no zero column.
\end{lemma}
\subsection{Paradoxical condition}
\begin{definition}\em
Let $AX=0$ be the homogenous system of equations corresponding to $Eq(\frak g,\mathcal E)$ and $B$ be as in  Lemma \ref{ll} and $L^j_i$ be the coefficient vector of the
equation
$$\sum_{x_j(C)\subseteq E_i}f_C =\sum_{x_0(C)\subseteq E_i}f_C = 0.$$
We say that $Eq\ge$ satisfies the {\it paradoxical condition} if each row of $B$ is of the
form
$\sum_{i=1}^m {R_i}$, where
$R_i\in \{L^j_{i},-L^j_{i}, L^j_i-L^k_i:\ \ \  1\leq j, k\leq n\}$ and
\\ \\
$$A=\left(
\begin{array}{c}
L_1^1 \\
\vdots \\
L^n_1 \\
\vdots\\
L^1_m \\
\vdots \\
L^n_m
\end{array}  \right)$$\\
In other words, each row of $B$ is of the form $\sum^m_{i=1}(L^{j_i}_i-L^{k_i}_i),$ for some $j_i, k_i\in\{0,\dots, n\}.$
\end{definition}

\begin{theorem} \label{pd} Let $G$ be a finitely generated group. Let $\frak g$ be an
ordered finite generating set for $G$ and $\mathcal E$ a finite partition of $G$. Suppose that the
associated system of configuration equations $Eq\ge$ admits no nonzero nonnegative
solution and satisfies the paradoxical condition. Then $G$ admits a paradoxical
decomposition in terms of $\frak g$ and $\mathcal E$.
\end{theorem}
A glance at the proof of the pervious theorem in \cite{y-r}, help us get a new  upper bound for the Tarski numbers of the groups under consideration. It is noticeable that in comparison with former results, this bound is not the most precise one, but is stated in terms of configurations and from this point of view is valuable to us. See the next theorem.

\begin{theorem} \cite[Corollary 3.7]{y-r} Let $\c$ be a configuration set of group $G$ satisfying the paradoxical
condition and $|\c|=l$.  Then  $t(G)\leq l + l^{2l}.$
\end{theorem}

\begin{remark}\em
In \cite{ytr}, the authors clarified the proof of the the theorem \ref{pd} and initiated the concept of {\it configuration graph}
to construct the paradoxical decomposition of a non-amenable discrete group using a system of configuration equations admitting no nonnegative non-zero solution. A short look at the proof of Theorem \ref{pd} shows that it is too complicated to work with. The configuration graph explains  this proof step by step and makes it applicable.
\end{remark}

\subsection{Normal condition}
Let $\pi:\{1,\dots,n\}\rightarrow \{1,\dots,n\}$ be a permutation
for the set $\{1,\dots,n\}$. Then $$P_{\pi}=\left(
\begin{array}{c}
e_{\pi(1)} \\
e_{\pi(2)} \\
\vdots \\
e_{\pi(n)}
\end{array} \right)$$ is
called the permutation matrix associated to $\pi$, where $e_{i}$
denotes the row vector of length $n$ with 1 in the i-th position
and 0 otherwise. 
 In the next definition $T$ is the following matrix
 
  $$T=\left(
\begin{array}{ccccc}
1& 0& 0& \dots & 0 \\
1& 1& 0& \dots& 0 \\
\vdots \\
1& 1& 1& \dots& 1
\end{array} \right).$$

\begin{definition}\em
Let $n\in \mathbb{N}$ and $\left( \begin{array}{cccc}
A_1 \\
A_2 \\
\vdots \\
A_n
\end{array} \right)
$ and $  \left( \begin{array}{cccc}
B_1 \\
B_2 \\
\vdots \\
B_n
\end{array} \right)
$ be two $(0,1)$-matrices with rows $A_i, B_i$. Let also the
vector $\sum_{i=1}^n (B_i-A_i)$ has strictly positive entries. If
there exists a permutation matrix $P_{\pi}$ such that the matrix
\begin{equation*}
TP_{\pi}(B-A)-P_{\pi(1\ 2\ \dots \ n)}A
\end{equation*}
has no entry less than $-1$, we say that the
homogenous system of equations $(B-A)X=0$ is {\it{normal.}}
\end{definition}
\begin{theorem}\label{main}
 If a subsystem of $Eq(\frak g,\mathcal E)$ is normal, then $G$
admits a paradoxical decomposition which is 
written in terms of $\frak g$ and $\mathcal E$.
\end{theorem}
For $Eq(\frak g,\mathcal E)$ satisfying the normal condition, there is a diagram that can give an upper bound for $\tau(G)$. The advantage of the second approach is to make the bound insofar as our choice of $\frak g$ and $\mathcal E$ permits (see \cite{arx}).

\section{Other algebraic structures}
Configurations are definable for some algebraic structures other than groups.
Abdollahi and Rejali introduced this concept for semigroups \cite{ar}.
They initiated paradoxical decompositions, generalized Theorem \ref{main} and extended the Tarski alternative for semigroups.
The notion of configuration for hypergroups can be found in \cite{bw}. In the latter paper for a hypergroup with a left Haar measure,
 the author gives the definition of configuration. Extending the last main results, he also provides a new characterization for left
 amenability of Lau algebras.\\

\subsection{Semigroups}
Throughout this subsection  $S$ denotes  a discrete semigroup. Recall that  $S$ is  left amenable, if there exists a left invariant mean
$M$ on $\ell^{\infty}(S)$ such that $M(xf)=M(f)$, for all $f\in \ell^{\infty}(S)$ and $x\in S$.
\begin{definition}\em
 Let  $\frak g=(g_1,\dots, g_n)$ be a generating sequence of elements of $S$ and
$\mathcal E=\{E_1,\dots,E_m\}$ be a partition of $S$. An $(n+1)$-tuple $C=(C_0,C_1,\dots,C_n),$
where $C_i\in \{1,\dots,m\}$ for each $i\in \{1,\dots,n\},$ is called a left configuration
corresponding to the configuration pair $\ge$, if there exists an element $x\in S$
with $x \in E_{C_0}$ such that

 $$g_ix\in E_{C_i}, \ \text{for\ each\ } i\in\{1, 2,\dots, n\}.$$
\end{definition}
 The set of all left
configurations corresponding to the configuration pair $\ge$ is denoted
by $Con_{\ell}\ge$.\\
Define $$x_0(C)=E_{C_0}\cap g_1^{-1} E_{C_1}\cap\dots\cap g_n^{-1} E_{C_n}$$
and $x_j(C) = g_jx_0(C)$, where for  $A \subseteq S$ and for a member $s$ in $S$ we put $$s^{-1}A :=\{t\in S: st\in A\}.$$

Completely similar to the group case, the left configuration equation corresponding to the configuration pair $\ge$ is
the system of equations

\begin{equation*}
\sum_{C\in Con_{\ell}(S)}\{f_C:\ x_j(C)\subseteq E_i\}=\sum_{C\in Con_{\ell}(S)}\{f_C:\ x_0(C)\subseteq E_i\},
\end{equation*}
where $i\in\{1,2,\dots,m\}$ and $j\in \{1,2,\dots,n\}.$ This system of equations is
denoted by $Eq_{\ell}\ge $.
A solution $(f_C)_{C\in Con_{\ell}\ge\}}$ to $Eq_{\ell}\ge$ satisfying $f_C\geq 0$, for\ all $C\in Con_{\ell}\ge$
and 
$$\sum_{C\in Con_{\ell}\ge} f_C = 1$$
is called
a normalized solution.\\

 In the following proposition, Theorem \ref{main}
is  generalized  for semigroups.\\

\begin{proposition} Let $S$ be a finitely generated discrete semigroup.The following statements are equivalent.\\
(i) $S$ is left amenable,\\
(ii) Each left configuration equation $Eq_{\ell}\ge$ has a normalised solution.
\end{proposition}

\begin{remark}\em
Paradoxical decompositions for semigroups can be defined and investigated out of the area of configurations. But because  of the similarity and relation between this concept and configurations in the case of groups, they are stated here.

Let $\{A_1,\dots,A_n;B_1,\dots,B_m\}$ be a partition of $S$ such that there exist
two subsets $\{g_1,\dots,g_n\}$ and $\{h_1,\dots,h_m\}$ of $S$ with the following property
\begin{eqnarray*}
S = \bigcup_{i=1}^ng_i^{-1} A_i=\bigcup_{j=1}^mh_j^{-1}B_j.
\end{eqnarray*}
Then we say that $S$ has a left paradoxical decomposition. \\

Using a method as in \cite[p.117]{paterson}, it is proved that a discrete semigroup  $S$ is left amenable if and only if there exists no left paradoxical decomposition for $S$. By this theorem which is a result in \cite{ar}, the authors extend the Tarski's theorem for semigroups which was asked in \cite[p. 120]{adler}. 
\end{remark}
\subsection{Hypergroups}
Let $H$ be a locally compact Hausdorff space and $M(H)$ be the space of finite regular Borel measures on $H$. The space $H$ is a hypergroup if there exists an associative binary
operation $*$ called convolution on $M(H)$ under which $M(H)$ is
an algebra. Moreover,
\begin{itemize}
\item for every $x, y$ in $H,$ $\delta_x*\delta_y$ is a
probability measure with
compact support.
  \item The mapping $(x,y)\mapsto \delta_x*\delta_y$ is a continuous map from $H\times H$ into
 $M(H)$ equipped with the weak* topology.
\item The mapping $(x,y)\mapsto supp(\delta_x*\delta_y)$ is a continuous mapping
from $H\times H$ into the compact subsets of $H$ equipped with the
Michael topology.
\item There exists a unique element  $e$ in $H$ such that
$\delta_e*\delta_x = \delta_x*\delta_e = \delta_x$ for all $x$ in
$H$.
\item There exists a homeomorphism $x\mapsto \check{x}$ of $H$
called involution satisfying $\check{\check x}= x$ for all $x\in
H$ and ${(\delta_x*\delta_y)\check{}} = \delta_{\check
y}*\delta_{\check x}$ for all $x,y\in H$, where
$\check{\mu}(A)=\mu(\check A),$ for any Borel subset $A$.
\item $e$ belongs to $supp(x*y)$ if and only if $y=\check{x}$,
\end{itemize}
where $\delta_x$ is the point mass measure at $x$. \\

For a measure $\mu\in M(H)$ and a Borel function $f$ on $H$, the left translation $\mu*f$ is defined by $\mu*f(x)=\check{\mu}*\delta_x(f)$. Recall that $H$ is amenable if there exist a normalized positive linear functional on $C(H)$, which is invariant under left translation. A Haar measure for $H$ is a non-zero positive regular Borel measure $\lambda$ which is left invariant; that is  $\lambda(\delta_x*f)=\lambda(f)$ for any $f\in C_C(H)$ and any $x\in H$.
\begin{definition}\cite{bw}  \em Let $H$ be a hypergroup with left Haar measure $\lambda$ and $\mathcal E=\{E_1,\dots,E_m\}$ be a finite measurable partition of $H$. Choose an
$n$-tuple of elements of $H$, $\frak h =\{h_1,\dots,h_n\}$. A
configuration is an $(n+1)$-tuple $C=(C_0,C_1,\dots,C_n)$ where
each $C_j\in \{1,\dots,m\}.$
\end{definition}

Despite of the group and semigroup cases,  configurations corresponding to $\mathcal E$ and $\frak h$ does only depend on the numbers $n$ and $m$, not on the elements of  $\mathcal E$ and $\frak h$.\\

Put $h_0 = e$ and for a fixed configuration $C$ define the positive 
real-valued function $\xi_0(C)$ on $H$ by
 $$\xi_0(C)(x) := \prod_{j=0}^n \delta_{h_j}*\delta_x(E_{C_j}).$$
 
\begin{definition} \em Let $\mathcal E$ and $\frak h$ be fixed and $\{z_C\ : C\in \c\}$ be variables
corresponding to the $m^{n+1}$ configurations. Consider the
$m\times n$ configuration equations $$\sum_{C\in \c, C_0=i} z_C =
\sum_{C\in \c, C_j=i} z_C$$ for each $i\in \{1,\dots,m\}$ and
$j\in\{1,\dots, n\}$. We say that a solution to these
configuration equations is positive if, for each $C\in \c,\
z_C\geq 0$;  normalized if $\sum_{C\in \c} z_C = 1$; and
inequality preserving if for every choice of $m^{n+1}$ real
numbers $\{\alpha_C :\ C\in \c\},$ $$ 0\leq \sum_{C\in \c}
\alpha_C\xi_0(C)\ a.e.\Rightarrow 0 \leq \sum_{C\in \c} \alpha_C\
z_C.$$
\end{definition}
Let the hypergroup $H$ be amenable with a left invariant mean $M$. suppose that $\frak h$ and $\mathcal E$ are given as before. For any configuration corresponding to  $(\frak h, \mathcal E),$ put $z_C=M(\xi_0(C)).$ Then we can see that $\{z_C,\ \ \ C \text{\ is\ a\ configuration}
 \}$ is a positive, normalized inequality preserving solution for the corresponding configuration equations.
The converse is also true but not so easy to reach. 

\begin{theorem}\cite[Theorm 2.6]{bw}\label{willson}
Let $H$ be a hypergroup with left Haar measure $\lambda$. $H$ is
amenable if and only if for all choices of $m, n, \mathfrak h$ and
$\mathcal E$ the $m\times n$ configuration equations have a
positive, normalized, inequality preserving solution.
\end{theorem}

\begin{remark}\em
It is not easy to extend the concept of paradoxical decomposition for hypergroups; But there is a special partition of unity for hypergroups which is called {\it paradoxical partition of unity}.
It is proved that a hypergroup is amenable if and only if it admits no paradoxical partition of unity (see \cite{yhypergroup}).
\end{remark}
\subsection{Configuration of Lau algebras}

Let $A$ be a Banach algebra and $X$ be a Banach $A$-bimodule. Then $X^*$ is a
Banach $A$-bimodule under the natural actions
$$\langle a.m, f \rangle=\langle m, f.a\rangle, \ \ \ \langle m.a,f \rangle=\langle m, a.f \rangle \ \ \ \ \ \ \ (a \in A, f\in X, m\in X^*).$$
A derivation $D: A\rightarrow X$ is a (bounded) linear map such that
$$D(ab)=D(a).b+a.D(b)\ \ \  (a, b \in A).$$
The derivation $D$ is inner if it is of the form $a\rightarrow a.f-f.a$ a for some $f\in X.$\\

A Lau algebra is a pair $(A,M)$ such that $A$ is a complex Banach algebra and $M$ is a $W^*$-algebra which $A$ is its predual, that is $A=M_*$
 and furthermore, the identity element of $M$, denoted by $e$, is a multiplicative linear functional on $A$. When there is no chance of confusion, we denote the pair $(A,M)$ simply by $A$. 
 Theses algebras were introduced by  Lau in 1983.
 In the original paper Lau called them $F$-algebras \cite{lau algebra}. This class of Banach algebras including group algebra, Fourier and Fourier
  Stieltjes algebra of a locally compact topological group and also the measure algebra of a locally compact hyperhroup has been of interest of
  many harmonic analysis researchers. In \cite{lau algebra} the left amenability of $F$-algebras is given and several characterization
 for left amenability is presented.
     Another characterization is based on the concept of configuration in a special sense which we shall review here.\\

The Lau algebra $A$ is called left amenable if for any two-sided Banach $A$-module $X$ such that
$$\phi.x=\phi(e)x\ \ \ \ \phi\in A, x\in X,$$
every bounded derivation from $A$ into $X^*$ is inner.\\

By $P(A)$ we denote the cone of all positive linear functionals in $A$ (as a subset of $A^{**}$). $P_1(A)$ denotes the functionals $\phi$ in $P(A)$ such that $\phi(e)=1$.
\begin{definition} \em
Let $(A,M)$ be a Lau algebra. Let $(\phi_1,\dots,\phi_n)\in (P_1(A))^n$ and $\{f_1,\dots,f_m\}\subseteq M$ such that each $f_i \geq 0$ and
$\sum_{i=1}^m f_i =e.$ We define an $(A,M)$-configuration as an ordered choice $C=(C_0,C_1,\dots,C_n)$
 with each $C_j\in\{1,\dots,m\}$ and define $\xi_0(C)$ by 
$$\xi_0(C)=\prod_{j=1}^n f_{C_j}.\phi_j.$$
\end{definition}
\noindent When the multiplication in $M$ is non-commutative, we
assume that the multiplication is done left to right as $j$ goes from 0 to $n$.\\

The configuration equations are defined 
in  variables $(z_C)_{C}$ by
$$\sum_{C_0=i}z_C =\sum_{C_j=i}z_C,\ \ \ (1\leq i\leq m, 1\leq j\leq n).$$
A solution to the configuration equations is again said to be positive if
each $z_C\geq 0$, normalized if
$\sum_C z_C = 1$ and inequality preserving if for
any choice of real numbers $\{\alpha_C\}$
$$0 \leq \sum_C\alpha_C \xi_0(C)\Rightarrow 0 \leq \sum_C\alpha_Cz_C.$$

\begin{theorem} A Lau algebra $(A,M)$ is left amenable if and only if
for all choices of $(\phi_1,\dots,\phi_n)\in(P_1(A))^n$ and $\{f_1,\dots,f_m\}\subseteq M$ such
that each $f_i\geq 0$ and
$\sum_{i=1}^mf_i= e$ the associated $(A,M)$-configuration
equations have a positive, normalized, inequality preserving solution.
\end{theorem}

\end{document}